\documentclass[11pt]{article}

\usepackage[]{graphicx,float,latexsym,times}
\usepackage{amsfonts,amstext,amsmath,amssymb,amsthm,bm}
\usepackage{caption2}

\long\def\symbolfootnote[#1]#2{\begingroup%
\def\thefootnote{\fnsymbol{footnote}}\footnote[#1]{#2}\endgroup}

\usepackage{titlesec} 

\titleformat{\section}{\large\bfseries}{\thesection.}{.5em}{}
\titlespacing*{\section}{0pt}{*3}{*2}
\titleformat{\subsection}{\normalfont\bfseries}{\thesubsection.}{.5em}{}
\titlespacing*{\subsection} {0pt}{*3}{*2}
\titleformat{\subsubsection}{\normalfont\bfseries}{\thesubsubsection.}{.5em}{}
\titlespacing*{\subsubsection} {0pt}{*3}{*2}

\theoremstyle{plain} 
\newtheorem{theorem}{Theorem}[section]
\newtheorem{lemma}{Lemma}[section]

\theoremstyle{definition} 

\newcommand{\bea}{\begin{eqnarray}}
\newcommand{\ena}{\end{eqnarray}}
\newcommand{\beas}{\begin{eqnarray*}}
\newcommand{\enas}{\end{eqnarray*}}
\newcommand{\beq}{\begin{equation}}
\newcommand{\enq}{\end{equation}}
\def\qed{\hfill \mbox{\rule{0.5em}{0.5em}}}

\newcommand{\ignore}[1]{}

\newcommand{\abs}[1]{\left|#1\right|}
\newcommand{\To}{\rightarrow}
\newcommand{\wtilde}[1]{\widetilde{#1}}
\newcommand{\qmq}[1]{\quad\mbox{#1}\quad}
\newcommand{\figformat}{pdf}

\numberwithin{equation}{section} 

\usepackage[authoryear]{natbib}

\begin{document}

\title{\textbf{\Large The Spend-It-All Region and Small Time Results for the Continuous Bomber Problem}}

\date{}

\maketitle


\author{
\begin{center}
\vskip -1cm

\textbf{\large Jay Bartroff$^1$, Larry Goldstein$^1$, and Ester Samuel-Cahn$^2$}

$^1$Department of Mathematics, University of Southern California, Los Angeles, California, USA\\
$^2$Department of Statistics and Center for the Study of Rationality,\\The Hebrew University of Jerusalem, Israel

\end{center}
}

\symbolfootnote[0]{\normalsize Address correspondence to Jay Bartroff,
Department of Mathematics, University of Southern California, KAP 108, Los Angeles,
CA 90089, USA; E-mail: bartroff@usc.edu}

{\small \noindent\textbf{Abstract:} A problem of optimally allocating partially effective
ammunition $x$ to be used on randomly arriving enemies in order to maximize an
aircraft's probability of surviving for time~$t$, known as the Bomber Problem, was
first posed by \citet{Klinger68}.  They conjectured a set of
apparently obvious monotonicity properties of the optimal allocation function $K(x,t)$.
Although some of these conjectures, and versions thereof, have been proved or disproved by other authors since
then, the remaining central question, that $K(x,t)$ is nondecreasing in~$x$, remains unsettled.
After reviewing the problem and summarizing the state of these conjectures, in the setting where $x$ is continuous we prove the existence of a ``spend-it-all'' region in which $K(x,t)=x$ and find its boundary, inside of which the long-standing, unproven conjecture of monotonicity of~$K(\cdot,t)$ holds. A new approach is then taken of directly estimating~$K(x,t)$ for small~$t$, providing a complete small-$t$ asymptotic description of~$K(x,t)$ and the optimal probability of survival.}

\section{INTRODUCTION}
\citet{Klinger68} introduced a problem of optimally allocating partially effective
ammunition to be used on enemies arriving at a Poisson rate in order to maximize the probability that an aircraft (hereafter ``the bomber'') survives for time~$t$, known as the Bomber Problem. Given an amount~$x$ of ammunition, let $K(x,t)$ denote the optimal amount of ammunition the bomber would use upon confronting an enemy at \textit{time}~$t$, defined as the time remaining to survive. The appearance of enemies is driven by a time-homogeneous Poisson process of known rate, taken to be $1$. An enemy survives the bomber's expenditure of an \hbox{amount~$y \in [0,x]$} of its ammunition with the geometric probability $q^y$, for some known $q\in(0,1)$, after which the enemy has a chance to destroy the bomber, which happens with known probability $v\in(0,1]$ (the $v=0$ case being trivial).  By rescaling $x$, we assume without loss of generality that $q=e^{-1}$, and hence the probability that the bomber survives an enemy encounter in which it spends an amount $y$ of its ammunition is
\begin{equation}\label{a-function}
a(y)=1-ve^{-y}.
\end{equation}
\citet{Klinger68} posed two seemingly obvious conjectures about the optimal allocation function $K(x,t)$:
\begin{enumerate}
\item[A:] $K(x,t)$ is nonincreasing in $t$ for all fixed $x\ge 0$;
\item[B:] $K(x,t)$ is nondecreasing in $x$ for all fixed $t\ge 0$.
\end{enumerate}
\citet{Klinger68} showed that [B] implies [A] when $v=1$, although, as will be discussed below, [B] remains in doubt. Improving the situation, \citet{Samuel70} showed that [A] holds without assuming [B] in the setting where units of ammunition $x$ are discrete, and in this setting also showed that a third conjecture holds:
\begin{enumerate}
\item[C:] The amount $x-K(x,t)$ held back by the bomber is nondecreasing in $x$ for all fixed $t\ge 0$.
\end{enumerate} [C] was first stated as a formal property by \citet{Simons90}, who claimed that it can be shown to hold for continuous $x$ and $t$ by arguments similar to the ones they provide for a case where both $x$ and $t$ are discrete, and they also make theoretical and computational progress toward~[B] in various discrete/continuous settings. Also in the setting where both $x$ and $t$ are continuous, Bartroff, Goldstein, Rinott, and Samuel-Cahn \citeyearpar{Bartroff09c} recently showed that [A] holds, and provide a full proof of [C] in this setting. \citet{Weber85} considered an infinite-horizon variant of the Bomber Problem in which the objective is to maximize the number of enemies shot down (thus removing $t$ from the problem) and found that, for discrete $x$, the property related to [B], that of monotonicity of $K(x)$, fails to hold. \citet{Shepp91} considered the infinite-horizon problem for continuous $x$ and reached the same conclusion.  On the other hand, \citet{Bartroff09c} consider the variation of the problem where the bomber is invincible, and both $x$ and $t$ are present and continuous, and show that [B] holds.

In spite of the results of \citet{Weber85}, \citet{Shepp91}, and \citet{Bartroff09c}, conjecture~[B] has not been settled in any close relative to the original Bomber Problem, and it remains the conjecture about which the least is known. To gain insight into the function $K(x,t)$, perhaps as a step towards resolving [B] in greater generality, we take a new approach to the Bomber Problem of directly estimating, or when possible solving for, $K(x,t)$ when both $x$ and $t$ are continuous. One might expect \hbox{\textit{a priori}} that if $x$ or $t$ is sufficiently small then the optimal strategy is to spend all or nearly all of the available ammunition~$x$, i.e., that $K(x,t)$ is equal to or nearly $x$. On the other hand, since the ammunition is assumed to be continuous it is not obvious that there exists a ``spend-it-all'' region where $K(x,t)$ is \textit{identically} $x$. In Section~\ref{sec:R1} we show that there is indeed a spend-it-all region of $(x,t)$ values for which $K(x,t)=x$ and where [B] holds, and we estimate the region's boundary in Theorem~\ref{thm:SIA}, and are able to find it exactly in most cases. However, in Section \ref{sec:smallt} we show that there are many other regimes in which $K(x,t)$ is not so simple, but can nevertheless be described asymptotically for small values of $t$. In particular, in Theorem~\ref{thm:smallt} we characterize the asymptotic behavior of $K(x,t)$ for small $t$ and show that regardless of how small $t$ is, there are large intervals of $x$ values for which $K(x,t)/x$ approaches any, even arbitrarily small, positive fraction, in stark contrast to the spend-it-all strategy. The relation of these results to the outstanding conjecture [B] and extensions are discussed in Section~\ref{sec:disc}.

\section{THE SPEND-IT-ALL REGION}\label{sec:R1}

In this section we describe an $(x,t)$-region where $K(x,t)$ is identically $x$, the so-called ``spend-it-all'' region. The boundary of this region is solved for, exactly as~(\ref{eq:fut}), except for a special configuration of the parameters $x, t, v$ in which the boundary is estimated from both sides; see~(\ref{eq:u01/2}). \citet{Bartroff10} has in the meantime shown that (\ref{eq:fut}) always gives the exact 
boundary of the spend-it-all region.

In what follows, let $u=1-v\in[0,1)$ denote the probability that the bomber survives an enemy's counterattack, let $P(x,t)$ denote the optimal probability of survival at time~$t$ when the bomber has ammunition~$x$, and let $H(x,t)$ denote the optimal conditional probability of survival given an enemy at time~$t$, with ammunition~$x$.

\bigskip

\begin{theorem}\label{thm:SIA} For $u\in(0,1)$ and $t>0$ define 
\begin{equation}\label{eq:fut}
f_u(t)=\log[1+u/(e^{tu}-1)],
\end{equation}
and extend this definition to $u=0$ by defining
$$f_0(t)=\lim_{u\To 0}f_u(t)=\log(1+t^{-1}).$$ For $u\in[0,1)$ and $t>0$ define 
\begin{equation}\label{eq:gut}
g_u(t)=\log(1+t^{-1}-u).
\end{equation}
If $u\in[0,1)$ and $t>0$ satisfy one of the following:
\begin{eqnarray}
&(i)&u=0,\label{eq:i}\\
&(ii)&u\in(0,1/2)\qmq{and} t\ge u^{-1}\log(2v),\label{eq:ii}\\
&(iii)&u\in[1/2,1),\label{eq:iii}
\end{eqnarray}
then 
\begin{equation}\label{eq:Kxiff}
\mbox{$K(x,t)=x$ if and only if $x\le f_u(t)$.}
\end{equation}
In the remaining case, where 
\begin{equation}\label{eq:iv}
u\in(0,1/2)\qmq{and} t<u^{-1}\log(2v),\end{equation}
we have 
\begin{equation}\label{eq:u01/2}
\mbox{$K(x,t)=x$ if $x\le g_u(t)$, and $K(x,t)<x$ if $x> f_u(t)$.}\end{equation}
\end{theorem}

\bigskip

The theorem may be summarized by saying that, except for the configuration of $t,u$ values in~(\ref{eq:iv}), the spend-it-all region's boundary is given exactly by $f_u(t)$, which is positive for all $t>0$ and approaches $0$ as $t\To\infty$. \citet{Bartroff10} has recently shown that $f_u(t)$ is the boundary of the spend-it-all region for all $t>0$ and $u\in[0,1)$.  Here, in the remaining case~(\ref{eq:iv}), the boundary is estimated from above by~$f_u(t)$ and from below by~$g_u(t)$, which is strictly less than $f_u(t)$ for all $t>0$ but asymptotically equivalent to it as $t\To 0$. Although $g_u(t)$ is negative for $t>u^{-1}$, it is utilized as a bound only when~(\ref{eq:iv}) holds, in which case $u^{-1}>u^{-1}\log(2v)>0$. A consequence of the theorem is that, regardless of the value of $u$, for any $x>0$ there is $t$ sufficiently small such that the optimal strategy spends it all (i.e., $K(x,t)=x$), and for any $t>0$ there is $x$ sufficiently small such that the optimal strategy spends it all.

\bigskip

\noindent\textit{Proof.}  We first prove that $K(x,t)=x$ when $x$ is bounded from above by $f_u(t)$ and one of (\ref{eq:i})-(\ref{eq:iii}) holds, or when~$x$ is bounded from above by $g_u(t)$ and~(\ref{eq:iv}) holds. To begin, fix $x,t$ and let $u$ be any value in $[0,1)$. We make use of the crude upper bound on the optimal survival probability
\begin{equation}\label{eq:Pu}
P(x,t)\le \exp(-vte^{-x})\qmq{for all} x,t>0,
\end{equation} which corresponds to the infeasible strategy of firing an amount $x$ of ammunition at every possible enemy, giving $$P(x,t)\le\sum_{i=0}^\infty e^{-t}[t a(x)]^i/i!=e^{-t}e^{t a(x)}=e^{-t(1-a(x))}=\exp(-vte^{-x}).$$ Using~(\ref{eq:Pu}), the optimal conditional survival probability is then $$H(x,t)=a(K(x,t))P(x-K(x,t),t)\le F(x-K(x,t)),$$ where for fixed $x$ and $t$ we write $$F(y)=a(x-y)\exp(-vte^{-y}).$$ By Lemma~\ref{lem:2stg} below, $F$ is unimodal on $\mathbb{R}$ with maximum at
$$y^*=\log\left(-vt+\sqrt{v^2t^2+4te^{x}}\right)-\log 2,$$ which is not necessarily in $[0,x]$.  In fact, if $x\le g_u(t)$, then
\begin{eqnarray*}
y^*&\le&\log\left(-vt+\sqrt{v^2t^2+4te^{g_u(t)}}\right)-\log 2\\
&=&\log\left(-vt+\sqrt{v^2t^2+4t(1+t^{-1}-u)}\right)-\log 2\\
&=&\log\left(-vt+\sqrt{(vt+2)^2}\right)-\log 2\\
&=&0,
\end{eqnarray*} hence in this case $\max_{y\in[0,x]} F(y)=F(0)=a(x)e^{-tv}$. If it were that $K(x,t)<x$, then we would have 
\begin{equation}\label{eq:H<SIA}
H(x,t)\le F(x-K(x,t))<F(0)=a(x)e^{-tv},\end{equation}
a contradiction since the latter is the conditional survival probability of the spend-it-all strategy:
\begin{equation}\label{eq:SIAH}
a(x)\sum_{i=0}^\infty u^ie^{-t}t^i/i! = a(x)e^{-t}e^{tu} = a(x)e^{-tv}.
\end{equation}
Note that $e^{-tv}$ is the probability of not being killed in the enemy's thinned Poisson process with parameter~$v$. The argument leading to~(\ref{eq:H<SIA}) thus shows that $K(x,t)=x$ whenever $x\le g_u(t)$; in particular,  $K(x,t)=x$ when~(\ref{eq:iv}) holds, or when~(\ref{eq:i}) holds after noting that $g_0(t)=f_0(t)$. For the remaining cases~(\ref{eq:ii}) and (\ref{eq:iii}), we obtain a tighter bound. Fix $x,t$ and let $u\in(0,1)$. Letting $$G(y)=a(x-y)e^{-t}[1+e^{vy/u}(e^{tu}-1)],$$ we claim that
 \begin{equation}
\label{eq:1}H(x,t)\le G(x-K(x,t)).
\end{equation} To prove (\ref{eq:1}), first, a simple verification yields that for any nonnegative $b_1,\ldots,b_i$, 
\begin{equation}
\label{eq:diveq}
\prod_{j=1}^i a(b_j) \le a(y/i)^n \quad \mbox{when $\sum_{j=1}^i b_j=y$.}
\end{equation}
Hence,
 $H(x,t)\le \wtilde{G}(x-K(x,t))$, where
$$\wtilde{G}(y)=a(x-y)e^{-t}\left[1+\sum_{i=1}^\infty \frac{(ta(y/i))^i}{i!}\right],$$ as the right hand side is the probability of survival for the infeasible strategy where one is given the number $i$ of future encounters, and divides the remaining amount $x-K(x,t)$ of ammunition optimally among them, firing $(x-K(x,t))/i$ at each. Next, we claim that
\begin{equation}
a(y/i)^i\le u^ie^{vy/u}\quad\mbox{for all $y\in[0,x]$ and all $i\ge 1$,}\label{eq:3}
\end{equation} implying that $\wtilde{G}(y)\le G(y)$ for all $y\in[0,x]$, and hence~(\ref{eq:1}). Letting $\rho_i=[a(y/i)/u]^i$, (\ref{eq:3}) is true since $\lim_{i\To\infty} \rho_i=e^{vy/u}$ and $\rho_i$ is evidently a nondecreasing sequence:
\begin{eqnarray*}
u^i(\rho_i-\rho_{i-1})&=&a(y/i)^i-u a(y/(i-1))^{i-1}\\
&=&a(y/i)^i-a(0)a(y/(i-1))^{i-1}\\
&\ge&0,
\end{eqnarray*} this last by (\ref{eq:diveq}). We will show below that if~(\ref{eq:ii}) or (\ref{eq:iii}) holds, then $G(y)$ is uniquely maximized over $y\in[0,x]$ at $y=0$. Since $G(0)=a(x)e^{-tv}$, it then follows that $K(x,t)=x$, as above. To verify the maximum of $G$, we show that $G'(0)\le 0$ and $G''(y)<0$ for all $y\in(0,x]$.  We compute
\begin{eqnarray*}
e^t G'(y)&=&-\frac{v}{u}\{e^{-x}[ue^y+e^{y/u}(e^{tu}-1)]-e^{vy/u}(e^{tu}-1)\}\\
e^t G''(y)&=&-\frac{v}{u^2}\{e^{-x}[u^2e^y+e^{y/u}(e^{tu}-1)]-ve^{vy/u}(e^{tu}-1)\}.
\end{eqnarray*} If $x\le f_u(t)$, which is equivalent to $e^{-x}\ge (1+u/(e^{tu}-1))^{-1}$, then we have
\begin{eqnarray*}
-\left(\frac{u}{v}\right)e^t G'(0)&=&e^{-x}(u+(e^{tu}-1))-(e^{tu}-1)\\
&\ge&\left(1+\frac{u}{e^{tu}-1}\right)^{-1}(e^{tu}-v)-(e^{tu}-1)\\
&=&\left(\frac{e^{tu}-1}{e^{tu}-v}\right)(e^{tu}-v)-(e^{tu}-1)\\
&=&0,
\end{eqnarray*} hence $G'(0)\le 0$.
 Next,
\begin{eqnarray*}
-\left(\frac{u^2 e^{-vy/u}}{v}\right)e^t G''(y)&=&e^{-x}[u^2e^{(2u-1)y/u}+e^y(e^{tu}-1)]-v(e^{tu}-1)\\
&=&e^{-x}p(y)-v(e^{tu}-1),
\end{eqnarray*} where $p(y)=u^2e^{(2u-1)y/u}+e^y(e^{tu}-1)$. When $u\ge 1/2$ the function $p(y)$ is clearly increasing in $y$ so for $y>0$ and $x\le f_u(t)$,
\begin{eqnarray}
-\left(\frac{u^2 e^{-vy/u}}{v}\right)e^t G''(y)&>& e^{-x}p(0)-v(e^{tu}-1)\nonumber\\
&=&e^{-x}(u^2+e^{tu}-1)-v(e^{tu}-1)\nonumber\\
&\ge&\left(\frac{e^{tu}-1}{e^{tu}-v}\right)(u^2+e^{tu}-1)-v(e^{tu}-1)\nonumber\\
&=&\left(\frac{e^{tu}-1}{e^{tu}-v}\right)[u^2+e^{tu}-1-v(e^{tu}-v)]\nonumber\\
&=&\left(\frac{e^{tu}-1}{e^{tu}-v}\right)[u(e^{tu}-2v)]\nonumber\\
&\ge&0, \label{eq:conc}
\end{eqnarray} since $u\ge 1/2$ implies that $2v\le 1\le e^{tu}$. Finally, we show that when~(\ref{eq:ii}) holds, $p(y)$ is still increasing. First compute
\begin{eqnarray*}
p'(y)&=&u(2u-1)e^{(2u-1)y/u}+e^y(e^{tu}-1),\\
p''(y)&=&(2u-1)^2e^{(2u-1)y/u}+e^y(e^{tu}-1)\quad\quad >0,
\end{eqnarray*} and
\begin{eqnarray*}
p'(0)&=&u(2u-1)+(e^{tu}-1)\\
&\ge&u(2u-1)+(2v-1)\quad\mbox{(since $t\ge u^{-1}\log(2v)$)}\\
&=&2u^2-3u+1\\
&=&2(1-u)(1/2-u)\\
&>&0
\end{eqnarray*} since $u<1/2$. Thus, the steps leading to~(\ref{eq:conc}) hold in this case as well, completing the proof that $K(x,t)=x$ when (\ref{eq:i}), (\ref{eq:ii}), (\ref{eq:iii}), or (\ref{eq:iv}) holds.

To complete the proof of the theorem, we show that $K(x,t)<x$ when $x>f_u(t)$. To do this, we bound $H(x,t)$ from below by the conditional survival probability $\underline{H}(y)$ of the strategy that fires an amount $y\in[0,x]$ of ammunition at the present enemy, fires all remaining ammunition $x-y$ at the next enemy (if one is encountered), and hopes for the best thereafter.  First assume that $u\in(0,1)$ and fix $x,t$ satisfying $x>f_u(t)$.  Then
\begin{eqnarray*}
\underline{H}(y)&=&a(y)\left[e^{-t}+e^{-t}\sum_{i=1}^\infty\frac{t^i a(x-y)u^{i-1}}{i!}\right]\\
&=&a(y)\left[e^{-t}+e^{-t}\frac{a(x-y)}{u}(e^{tu}-1)\right]\\
&=&e^{-t}a(y)\left[1+\left(\frac{e^{tu}-1}{u}\right)a(x-y)\right].
\end{eqnarray*} By applying Lemma~\ref{lem:2stg} with $A=(e^{tu}-1)/u$,  we see that $\underline{H}(y)$ is unimodal with maximum at $K^*(x,t)=(x+f_u(t))/2$, which, since $x>f_u(t)$, satisfies $K^*(x,t)<(x+x)/2=x.$ If it were that $K(x,t)=x$, then we would have $$H(x,t)=a(x)e^{-tv}=\underline{H}(x)<\underline{H}(K^*(x,t)),$$ a contradiction. If $u=0$, the conditional survival probability of this strategy is
$$\underline{H}(y)=a(y)[e^{-t}+e^{-t}t a(x-y)]=e^{-t}a(y)[1+ta(x-y)],$$
and a similar argument applies: By Lemma~\ref{lem:2stg} with $A=t$, the function~$\underline{H}(y)$ is unimodal with maximum at $K^*(x,t)=(x+f_0(t))/2< x$, leading to the same contradiction. \qed

\bigskip

\begin{lemma}\label{lem:2stg} Fix $x>0$, $t>0$, and $v\in(0,1]$. The function 
\begin{equation}\label{eq:F1}
y\mapsto a(x-y)\exp(-vte^{-y})
\end{equation} is unimodal on $\mathbb{R}$ with maximum at \begin{equation}\label{eq:F1max}
y^*=\log\left(-vt+\sqrt{v^2t^2+4te^{x}}\right)-\log 2.
\end{equation} 
For any fixed $A>0$, the function 
\begin{equation}\label{eq:F2}y\mapsto a(y)[1+A a(x-y)]\end{equation} is unimodal on $\mathbb{R}$ with maximum at 
\begin{equation}\label{eq:F2max}y^*=[x+\log(1+A^{-1})]/2.\end{equation}
\end{lemma}

\noindent\textit{Proof.} Taking the derivative of~(\ref{eq:F1}) with respect to $y$ and setting $z=e^y$ gives 
\begin{equation}\label{eq:dF1}
-ve^{-x-y}\exp(-vte^{-y})(e^{2y}+vte^y-te^x)=-ve^{-x-y}\exp(-vte^{-y})(z^2+vtz-te^x).
\end{equation}
Since $z>0$, the function~(\ref{eq:F1}) increases in $y=\log z$ up to the log of the positive root of the quadratic in~(\ref{eq:dF1}), which is~(\ref{eq:F1max}), and decreases thereafter. Similarly, the derivative of~(\ref{eq:F2}) with respect to~$y$ is $$-v(Ae^{-x+y}-(1+A)e^{-y})=-ve^{-y}(Ae^{-x}z^2-(1+A)),$$ and solving for the root gives~(\ref{eq:F2max}).\qed

\section{AN ASYMPTOTIC CHARACTERIZATION OF $K(x,t)$}\label{sec:smallt}

In this section we give an asymptotic description of the optimal allocation function~$K(x,t)$ as $t\To 0$, and for this it suffices to consider sequences $(x,t)$ with $t\To 0$.  In addition to giving an asymptotic description of the optimal survival probability~$P(x,t)$ and the optimal conditional survival probability~$H(x,t)$, our main goal is to characterize the fraction $K(x,t)/x$ of the current ammunition $x$ spent by the optimal strategy at time~$t$, and it turns out that $K(x,t)/x$ approaches a finite nonzero limit on sequences $(x,t)$ such that $\abs{\log t}/x$ approaches a finite nonzero limit.  We thus give an essentially complete asymptotic description of $K(x,t)$ by considering sequences $(x,t)=(x_t,t)$ such that
\begin{equation}
\frac{\abs{\log t}}{x}\To\rho\in(0,\infty)\quad\mbox{as $t\To 0$,}\label{eq:logt/x}
\end{equation} leaving divergent sequences to be handled by considering subsequences. We will write $x=x_t$ when we wish to emphasize the dependence of~$x$ on~$t$, but most of the time this notation will be suppressed. Note that a consequence of~(\ref{eq:logt/x}) is that $x\To\infty$ at the same rate at which $\abs{\log t}\To\infty$ as $t\To 0$. It should perhaps not be surprising that this is the nontrivial asymptotic regime since the boundary of the spend-it-all region found in Theorem~\ref{thm:SIA} is asymptotically equivalent to~$\abs{\log t}$ as $t\To 0$. In what follows, let ${1\choose 2}^{-1}=\infty$.

\bigskip

\begin{theorem}\label{thm:smallt} Under~(\ref{eq:logt/x}), let $j\in\{1,2,\ldots\}$ be such that
\begin{equation}
\label{eq:Ij}
{j+1\choose 2}^{-1}\le \rho<{j\choose 2}^{-1}.\end{equation} Then, as $t\To 0$,
\begin{eqnarray}
\frac{K(x,t)}{x}&\To& 1/j+\rho(j-1)/2\label{eq:K/x}\\
\frac{1}{x}\abs{\log(1-H(x,t))}&\To& 1/j+\rho(j-1)/2\label{eq:H/x}\\
\frac{1}{x}\abs{\log(1-P(x,t))}&\To& 1/j+\rho(j+1)/2.\label{eq:P/x}
\end{eqnarray}
\end{theorem}

\bigskip

The theorem is proved in the next subsection. First, we briefly discuss the result. Note that the~$j$ satisfying~(\ref{eq:Ij}) is nonincreasing in $\rho$ and, in particular, $\rho\ge 1$ corresponds to $j=1$ while $\rho<1$ corresponds to $j>1$. The right hand sides of (\ref{eq:K/x}) and (\ref{eq:H/x}) equal 1 for $j=1$, and are in the interval $[2/(j+1),2/j)$ for $j\ge 2$; similarly, the right hand side of~(\ref{eq:P/x}) is in the interval $[2/j,2/(j-1))$ for all $j\ge 1$. In particular, (\ref{eq:K/x}) implies that $K(x,t)/x$ can take on any value in $(0,1]$. The rates of convergence in (\ref{eq:K/x})-(\ref{eq:P/x}) are functions of the rate of convergence in~(\ref{eq:logt/x}). Specifically, without assuming more than $\abs{\log t}-\rho x=o(x)$ in~(\ref{eq:logt/x}), the same $o(x)$ term appears in the convergence of $K(x,t)$,  $\abs{\log(1-H(x,t))}$, and $\abs{\log(1-P(x,t))}$ in (\ref{eq:K/x})-(\ref{eq:P/x}). However, when $\rho>1$, the convergence is $O(1/x)$ in (\ref{eq:K/x}) and (\ref{eq:H/x}), but in no other cases, an artifact of the natural upper bound $K(x,t)\le x$ that is relevant only in the $\rho>1$ case.

The result~(\ref{eq:K/x}) can equivalently be stated as, under~(\ref{eq:logt/x}),
\begin{equation}\label{eq:Ksim}
K(x,t)\sim \frac{x}{j}+\left(\frac{j-1}{2}\right)\abs{\log t}\end{equation} as $t\To 0$ for $j$ satisfying (\ref{eq:Ij}). Hence, for small $t$, the first quadrant of the $(x,t)$-plane can be thought of as partitioned into the regions
\begin{equation}\label{eq:Rj}
R_j=\left\{(x,t):\quad x>0,\quad t>0,\quad {j+1\choose 2}^{-1}\le \frac{\abs{\log t}}{x}<{j\choose 2}^{-1}\right\},\quad j=1,2,\ldots,
\end{equation} which determine the asymptotic behavior of the optimal strategy.  Figure~1 plots~(\ref{eq:Ksim}) and the boundaries of the first few $R_j$. Note that although~(\ref{eq:Ksim}) varies smoothly within each $R_j$, it is continuous but not smooth at the lower boundary of $R_j$. For small $t$, $K(x,t)$ given by~(\ref{eq:K/x}) turns out to be such that if $(x,t)\in R_j$, then after firing $K(x,t)$ at an immediate enemy, the new state $(x-K(x,t),t)$ lies in $R_{j-1}$. This leads to the inductive method of proof, given in the next section. The boundary of the $R_1$ region is asymptotically equivalent to the estimates of the spend-it-all region's boundary in Theorem~\ref{thm:SIA} in the strong sense that their difference is $o(1)$ as $t\To 0$.

\begin{figure}
\begin{center}\scalebox{.4}{\includegraphics{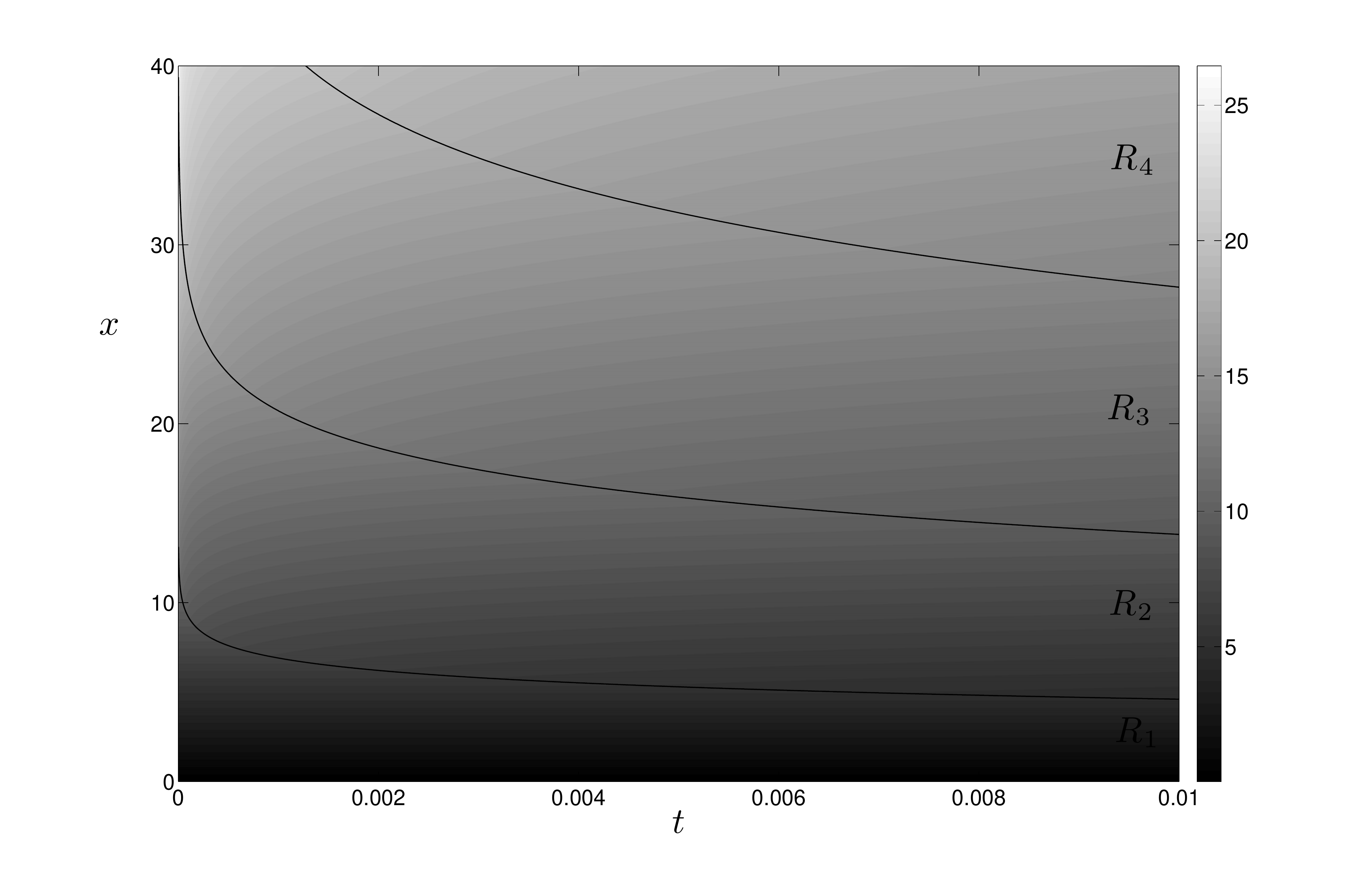}}\\
\textit{\textbf{Figure 1.}} The small-$t$ asymptotic approximation~(\ref{eq:Ksim}) of $K(x,t)$.
\end{center}
\end{figure}

\subsection{Proof of Theorem \ref{thm:smallt}}\label{sec:smalltpf}
We proceed by induction on $j$. To begin, assume that $j=1$. The optimal conditional survival probability $H(x,t)$ is bounded below by the conditional survival probability of the spend-it-all strategy~(\ref{eq:SIAH}), giving
\begin{eqnarray}
H(x,t)&\ge&a(x)e^{-tv}\nonumber\\
&\ge&(1-ve^{-x})(1-tv)\nonumber\\
&=&1-ve^{-x}-ve^{-\rho x+o(x)}+v^2 e^{-(\rho+1)x+o(x)}.\label{eq:j1Kl}
\end{eqnarray} 
Moving to $P(x,t)$, it is bounded below by the survival probability of the strategy
that fires~$x$ at the first enemy.  Under such a strategy, the bomber survives if no enemy plane
arrives during time~$t$, which happens with probability~$e^{-t}$, or if the bomber encounters and survives one enemy, which
happens with probability~$t e^{-t}a(x)$, and  ignoring other  enemy encounters we obtain
\begin{eqnarray}
P(x,t)&\ge&e^{-t}[1+ta(x)]\nonumber\\
&\ge&(1-t)[1+t-vte^{-x}]\nonumber\\
&=&1-vte^{-x}-t^2+vt^2e^{-x}\nonumber\\
&=&1-ve^{-(\rho+1)x+o(x)}-e^{-2\rho x+o(x)}+ve^{-(2\rho+1)x+o(x)}\nonumber\\
&\ge&1-ve^{-(\rho+1)x+o(x)}.\label{eq:j1Pl}
\end{eqnarray} On the other hand, $P(x,t)$ is bounded above by the survival probability of the infeasible strategy that fires $x$ at the first enemy and, upon survival of this encounter, is guaranteed survival thereafter, so
\begin{eqnarray}
P(x,t)&\le& e^{-t}[1+a(x)(e^t-1)]\nonumber\\
&=&1-ve^{-x}+ve^{-x}e^{-t}\nonumber\\
&\le&1-ve^{-x}+ve^{-x}(1-t+t^2/2)\nonumber\\
&=&1-vte^{-x}+vt^2e^{-x}/2\nonumber\\
&=&1-ve^{-(\rho+1)x+o(x)}+ve^{-(2\rho+1)x+o(x)}/2\nonumber\\
&=&1-e^{-(\rho+1)x+o(x)}.\label{eq:j1Pu}
\end{eqnarray}
For this $j=1$ case, we consider separately the cases $\rho=1$ and $\rho\in(1,\infty)$. First assume that $\rho>1$. In this case, (\ref{eq:j1Kl}) is $1-ve^{-x}(1+o(1))$, and as the conditional probability of the spend-it-all strategy is bounded above by $a(K(x,t))\cdot 1$, the probability of surviving the first encounter when expending the optimal amount $K(x,t)$ and ignoring future danger, we find that $1-ve^{-x}(1+o(1))\le 1-ve^{-K(x,t)}$ so that $x+o(1)\le K(x,t)\le x$. To estimate $H(x,t)$, plug $K(x,t)=x+o(1)$ into $a(K(x,t))$ to get $$a(K(x,t))=1-ve^{-x}(1+o(1)).$$ This is the same order as the lower bound~(\ref{eq:j1Kl}), hence $H(x,t)=1-ve^{-x}(1+o(1))$, which implies~(\ref{eq:H/x}) in this case. The limit~(\ref{eq:P/x}) holds as well in this case since both (\ref{eq:j1Pl}) and (\ref{eq:j1Pu}) are of order $$1-e^{-(\rho+1)x+o(x)}.$$ Since $H(x,t)=1-ve^{-x}(1+o(1))$ is equivalent to
$$\frac{1}{x}\abs{\log(1-H(x,t))}=1-(\log v)/x+o(1/x)=1+O(1/x),$$ the error term on the right hand side of~(\ref{eq:H/x}) in this case is $O(1/x)$; this holds for~(\ref{eq:K/x}) and (\ref{eq:H/x}) when $\rho>1$, but for no other cases.

Now let $\rho=1$. The lower bound~(\ref{eq:j1Kl}) is
\begin{equation}\label{eq:j11Kl}1-2ve^{-x+o(x)}=1-e^{-x+o(x)},\end{equation} and by Lemma~\ref{lem:Kl} below we have $x+o(x)\le K(x,t)\le x$. Plugging $K(x,t)=x+o(x)$ into the upper bound $a(K(x,t))$ gives $H(x,t)\le1-e^{-x+o(x)}$, the same order as the lower bound~(\ref{eq:j11Kl}), hence $$\frac{1}{x}\abs{\log(1-H(x,t))}\To 1.$$ The lower bound~(\ref{eq:j1Pl}) gives
$$P(x,t)\ge 1-e^{-2x+o(x)},$$ and the upper bound~(\ref{eq:j1Pu}) gives the same order, hence $$\frac{1}{x}\abs{\log(1-P(x,t))}\To 2=\rho+1.$$ This concludes the $j=1$ case.

Now let $I_j$ denote the half-closed interval~(\ref{eq:Ij}), i.e.,
\begin{equation}\label{eq:Ij2}
I_j=\left[{j+1\choose 2}^{-1}, {j\choose 2}^{-1}\right),
\end{equation}
and let $\alpha_j(\rho)$ and $\beta_j(\rho)$ denote the right hand sides of~(\ref{eq:K/x}) and (\ref{eq:P/x}), respectively, i.e.,
\begin{eqnarray}
\alpha_j(\rho)&=&1/j+\rho(j-1)/2\label{eq:aj}\\
\beta_j(\rho)&=&1/j+\rho(j+1)/2.\label{eq:bj}
\end{eqnarray}
For the inductive step, assume that (\ref{eq:K/x})-(\ref{eq:P/x}) hold for $j$ and let $\rho$ belong to $I_{j+1}$. $H(x,t)$ is bounded below by the conditional survival probability of the strategy $\underline{H}(x,t)$ that fires $\wtilde{K}(x)=\alpha_{j+1}(\rho)x$ at the first enemy, and then behaves optimally thereafter. Letting $$x':=x-\wtilde{K}(x)=x[1-\alpha_{j+1}(\rho)],$$ we have
$$\rho':=\lim_{t\To 0}\frac{\abs{\log t}}{x'}=\frac{\rho}{1-\alpha_{j+1}(\rho)}\in I_j$$ by Lemma \ref{lem:alphaj} below. Then, by the inductive hypothesis, we have
\begin{equation}\label{eq:j2Hl}
\underline{H}(x,t) = a(\wtilde{K}(x))P(x',t) = [1-ve^{-\alpha_{j+1}(\rho)x}][1-e^{-\beta_{j}(\rho')x'+o(x')}],
\end{equation} and
$$\beta_{j}(\rho')\frac{x'}{x}=\beta_{j}(\rho')[1-\alpha_{j+1}(\rho)]=\alpha_{j+1}(\rho)$$ by Lemma~\ref{lem:alphaj}, giving
\begin{equation}
\underline{H}(x,t)=[1-e^{-\alpha_{j+1}(\rho)x+o(x)}]^2=1-e^{-\alpha_{j+1}(\rho)x+o(x)}.\label{eq:j+1Hl}
\end{equation} Lemma \ref{lem:Kl} then implies that
\begin{equation}
\label{eq:j+1Kl}
K(x,t)\ge \alpha_{j+1}(\rho)x+o(x),\end{equation} and we will show that this expression actually holds with equality. To do this, we consider sequences $(x,t)$ still for which $\abs{\log t}/x\To \rho\in I_{j+1}$ and on which $$\tau:=\lim_{t\To 0} \frac{K(x,t)}{x}$$ exists, and we will show that $\tau=\alpha_{j+1}(\rho)$ is the only possible limit.  This suffices to show that the $\limsup$ and $\liminf$ of $K(x,t)/x$ both equal $\alpha_{j+1}(\rho)$.

By~(\ref{eq:j+1Kl}), we know that the only possible values of $\tau$ lie in $[\alpha_{j+1}(\rho),1]$. First, suppose that there is a sequence $(x,t)$ on which $\tau\in(\alpha_{j+1}(\rho),1)$. Then $$x'':=x-K(x,t)\sim(1-\tau)x\qmq{and} \rho'':=\lim_{t\To 0}\frac{\abs{\log t}}{x''}=\frac{\rho}{1-\tau}>\frac{\rho}{1-\alpha_{j+1}(\rho)}=\rho'\in I_j$$ by Lemma~\ref{lem:alphaj}, so let $i\in\{1,2,\ldots,j\}$ be such that $\rho''\in I_i$. Then, again by the inductive hypothesis, we would have
\begin{equation}\label{eq:j+1Hu3}
H(x,t)=a(K(x,t))P(x'',t)=[1-ve^{-\tau x+o(x)}][1-e^{-\beta_i(\rho'')x''+o(x'')}],
\end{equation} and
\begin{eqnarray}
\beta_i(\rho'')\frac{x''}{x}=\left[\frac{1}{i}+\frac{\rho''(i+1)}{2}\right]\frac{x''}{x}&\To& \left[\frac{1}{i}+\frac{\rho(i+1)}{2(1-\tau)}\right](1-\tau) \nonumber\\
&=& \frac{1-\tau}{i}+\frac{\rho(i+1)}{2} \label{eq:j+1Hu4}.
\end{eqnarray}
If $i<j$, then $\rho/(1-\tau)\in I_i$ implies that $(1-\tau)\le \rho{i+1\choose 2}$, so (\ref{eq:j+1Hu4}) becomes
\begin{eqnarray*}
\frac{1-\tau}{i}+\frac{\rho(i+1)}{2}&\le&\rho(i+1)\nonumber\\
&=&\rho(i+1-j/2)+\rho j/2\nonumber\\
&<&{j+1\choose 2}^{-1}(j-j/2)+\rho j/2\quad\mbox{(since $\rho\in I_{j+1}$ and $i< j$)}\nonumber\\
&=&\alpha_{j+1}(\rho).
\end{eqnarray*} If $i=j$, then~(\ref{eq:j+1Hu4}) becomes
\begin{eqnarray*}
\frac{1-\tau}{j}+\frac{\rho(j+1)}{2}&<&\frac{1-\alpha_{j+1}(\rho)}{j}+\frac{\rho(j+1)}{2}\\
&=&\alpha_{j+1}(\rho)-\frac{[(j+1)\alpha_{j+1}(\rho)-1]}{j}+\frac{\rho(j+1)}{2}\\
&=&\alpha_{j+1}(\rho).
\end{eqnarray*} In both cases we have shown that~(\ref{eq:j+1Hu4}) is less than $\alpha_{j+1}(\rho)<\tau$, which implies that (\ref{eq:j+1Hu3}) is $$1-\exp[-((1-\tau)/i+\rho(i+1)/2)x+o(x)]$$ and is hence smaller than~(\ref{eq:j+1Hl}) for small $t$, a contradiction.

Now assume that there is a sequence $(x,t)$ on which $\tau=1$. Using the crude bound $a(K(x,t))\le 1$ and~(\ref{eq:Pu}),
\begin{eqnarray*}
H(x,t)&=&a(K(x,t))P(x-K(x,t),t)\\
&\le&1\cdot \exp[-tve^{-(x-K(x,t))}]\\
&\le&1-vte^{-(x-K(x,t))}+v^2t^2e^{-2(x-K(x,t))}/2\\
&=&1-ve^{-\rho x+o(x)}+v^2e^{-2\rho x+o(x)}\\
&=&1-e^{-\rho x+o(x)},
\end{eqnarray*} which leads to the same contradiction since $\rho<\alpha_{j+1}(\rho)$ by Lemma~\ref{lem:alphaj}. We have shown that $\alpha_{j+1}(\rho)$ is the only possible value of $\tau$, hence~(\ref{eq:j+1Kl}) holds with equality and~(\ref{eq:j+1Hl}) holds for $H(x,t)$. All that remains is to verify~(\ref{eq:P/x}) for the $j+1$ case.

Let $T$ denote the exponentially distributed waiting time, with mean~$1$, until the first enemy, and recall that we write $x=x_t$ to emphasize the dependence on $t$.  Then $P$ and $H$ are related through the expectation
\begin{eqnarray}
P(x,t)=P(x_t,t)&=&E\left[H(x_t,t-T)\bm{1}\{T<t\}+\bm{1}\{T\ge t\}\right]\nonumber\\
&=&\int_0^t H(x_t,t-r)e^{-r} dr+P(T\ge t)\nonumber\\
&=&e^{-t}\left[\int_0^t H(x_t,s)e^s ds+1\right].\label{eq:Pint}
\end{eqnarray} Using~(\ref{eq:Pint}) and that $H(x,\cdot)$ is nonincreasing, we have
\begin{eqnarray*}
P(x,t)&\ge&e^{-t}\left[H(x_t,t)\int_0^t e^s ds+1\right]\\
&=&e^{-t}\left[H(x_t,t)(e^t-1)+1\right]\\
&=&e^{-t}[1-H(x_t,t)]+H(x_t,t)\\
&\ge&(1-t)[1-H(x_t,t)]+H(x_t,t)\\
&=&1-t[1-H(x_t,t)]\\
&=&1-e^{-\rho x+o(x)}e^{-\alpha_{j+1}(\rho)x+o(x)}\\
&=&1-e^{-[\rho+\alpha_{j+1}(\rho)]x+o(x)}\\
&=&1-e^{-\beta_{j+1}(\rho)x+o(x)}\\
\end{eqnarray*} by Lemma \ref{lem:alphaj}. We bound $P(x,t)$ from above by a function of the same order. Fix $\delta\in(0,1)$ and note that
\begin{equation}\label{eq:dt}
\frac{\abs{\log(\delta t)}}{x_t}=\frac{-\log(\delta t)}{x_t}=\frac{-\log t}{x_t}+\frac{- \log \delta}{x_t}=\frac{\abs{\log t}}{x_t}+\frac{\abs{\log \delta}}{x_t}\To \rho.
\end{equation} Then, by (\ref{eq:Pint}),
\begin{eqnarray*}
P(x,t)&\le&e^{-t}\left[\int_0^{\delta t} e^s ds+H(x_t,\delta t)\int_{\delta t}^t e^s ds+1\right]\\
&=&e^{-(1-\delta)t}[1-H(x_t,\delta t)]+H(x_t,\delta t)\\
&\le&[1-(1-\delta)t+t^2][1-H(x_t,\delta t)]+H(x_t,\delta t)\\
&=&1-(1-\delta)t[1-H(x_t,\delta t)]+t^2[1-H(x_t,\delta t)]\\
&=&1-(1-\delta)e^{-\rho x+o(x)}e^{-\alpha_{j+1}(\rho)x+o(x)}+e^{-2\rho x+o(x)}e^{-\alpha_{j+1}(\rho)x+o(x)}\quad\mbox{(by (\ref{eq:dt}))}\\
&=&1-e^{-[\rho+\alpha_{j+1}(\rho)]x+o(x)}\\
&=&1-e^{-\beta_{j+1}(\rho)x+o(x)},\\
\end{eqnarray*} completing the proof of Theorem~\ref{thm:smallt}, except for the following lemmas. The first collects various facts relating $\alpha_j(\rho)$, $\beta_j(\rho)$, and $\rho$, and the second provides a crude but useful bound on $K(x,t)$.

\bigskip

\begin{lemma}\label{lem:alphaj} Let $I_j$, $\alpha_j(\rho)$, and $\beta_j(\rho)$ be as in (\ref{eq:Ij2})-(\ref{eq:bj}). Assume that $\rho\in I_{j+1}$ for some $j\ge 1$, and let $\rho'=\rho/[1-\alpha_{j+1}(\rho)]$. Then
\begin{eqnarray}
\rho&<&\alpha_{j+1}(\rho),\label{eq:alphj1}\\
\rho'&\in& I_{j},\label{eq:alphj2}\\
\beta_{j}(\rho')&=&[1/\alpha_{j+1}(\rho)-1]^{-1}.\label{eq:alphj3}\\
\alpha_{j+1}(\rho)+\rho&=&\beta_{j+1}(\rho)\label{eq:alphj4}
\end{eqnarray}
\end{lemma}

\noindent\textit{Proof.} Let $\rho\in I_{j+1}$. Then
\begin{eqnarray*}
\beta_j(\rho')&=&\frac{1}{j}+\frac{\rho'(j+1)}{2}\\
&=&\frac{1}{j}+\frac{2(j+1)^2}{2j(2/\rho-(j+1))}\\
&=&\frac{\alpha_{j+1}(\rho)}{1-\alpha_{j+1}(\rho)}
\end{eqnarray*} after some simplifying, proving (\ref{eq:alphj3}). For~(\ref{eq:alphj1}),
\begin{eqnarray*}
\rho&=&\rho(1-j/2)+\rho j/2\\
&\le&\left\{
\begin{array}{ll}
 \rho/2+\rho/2, &j=1 \\
  \rho j/2,& j\ge 2
  \end{array}
\right.\\
&<&\left\{
\begin{array}{ll}
 1/2+\rho/2, &j=1 \\
 1/(j+1)+ \rho j/2,& j\ge 2
  \end{array}
\right.\\
&=&\alpha_{j+1}(\rho).
\end{eqnarray*}

For~(\ref{eq:alphj2}),
\begin{eqnarray*}
\rho'&=&\frac{2(j+1)}{j[2/\rho-(j+1)]}\\
&\in&\left[\frac{2(j+1)}{j[2{j+2\choose 2}-(j+1)]},\frac{2(j+1)}{j[2{j+1\choose 2}-(j+1)]}\right)\\
&=&\left[{j+1\choose 2}^{-1},{j\choose 2}^{-1} \right)=I_j.
\end{eqnarray*}
For (\ref{eq:alphj4}), $\alpha_{j+1}(\rho)+\rho=1/(j+1)+\rho(j+2)/2=\beta_{j+1}(\rho)$.
\qed

\bigskip

\begin{lemma}\label{lem:Kl} If there is a $\gamma\in(0,1]$ such that $H(x,t)\ge 1-e^{-\gamma x+o(x)}$, then $K(x,t)\ge\gamma x+o(x)$.
\end{lemma}

\noindent\textit{Proof.} We have $$H(x,t)=a(K(x,t))P(x-K(x,t),t)\le a(K(x,t))\cdot 1=1-ve^{-K(x,t)},$$ and setting this last $\ge$ the assumed lower bound $1-e^{-\gamma x+o(x)}$ leads to $K(x,t)\ge\gamma x+o(x)$.\qed

\section{DISCUSSION}\label{sec:disc}

In Section \ref{sec:smallt} an inductive method is used to estimate the limiting optimal fraction $K(x,t)/x$ of ammunition used as $t\To 0$.  The same result holds when the bomber is restricted to only firing discrete units (integers, say) of ammunition $x$, the only modification of the proof needed is to replace $x$ by $\lfloor x\rfloor$ (the largest integer $\le x$) in the appropriate places. For example, in the $\rho>1$ case in the proof of Theorem~\ref{thm:smallt}, we have $H(x,t)\ge a(\lfloor x\rfloor)e^{tv}$, which leads to $\lfloor x\rfloor+o(1)\le K(x,t)\le \lfloor x\rfloor$, and hence $K(x,t)/x\To 1$, using that $\lfloor x\rfloor/x\To 1$.

Theorem~\ref{thm:SIA} shows that $K(x,t)=x$ in a region asymptotically equivalent to $R_1$ in (\ref{eq:Rj}) and, this being monotone in $x$, that conjecture [B] holds in this region. It is therefore natural to ask if the estimates of $K(x,t)$ in $R_j$ given by Theorem~\ref{thm:smallt} can be used to shed any light on conjecture~[B] for $j\ge 2$. One thing we can say is that [B] is satisfied \textit{in the limit} as $t\To 0$ in the following sense. Letting $x_1\le x_2$ be such that $\lim_{t\To 0} \abs{\log t}/x_1\in R_j$ and $\lim_{t\To 0} \abs{\log t}/x_2\in R_{j'}$ for some $j\le j'$, by (\ref{eq:Ksim}) we have
\begin{equation}\label{eq:Kdiff}
K(x_2,t)-K(x_1,t)\sim\frac{x_2}{j'}+\left(\frac{j'-1}{2}\right)\abs{\log t}-\left[\frac{x_1}{j}+\left(\frac{j-1}{2}\right)\abs{\log t}\right].\end{equation}
If $j=j'$, then (\ref{eq:Kdiff}) is $(x_2-x_1)/j\ge 0$. If $j<j'$, then~(\ref{eq:Kdiff}) divided by $\abs{\log t}$ is
\begin{eqnarray*}
\frac{x_2}{j' \abs{\log t}}-\frac{x_1}{j \abs{\log t}}+\frac{j'-j}{2}&>& \left(\frac{{j'\choose 2}}{j'}-\frac{{j+1\choose 2}}{j }+\frac{j'-j}{2}\right)(1+o(1))\\
&=&(j'-j-1)(1+o(1))
\end{eqnarray*} which approaches a nonnegative limit. However, to make this argument hold for, say, all $x_1\le x_2$ sufficiently large and all $t$ sufficiently small, higher order asymptotics are needed. In particular, the rate of convergence in (\ref{eq:K/x}) as a function of $x$ and $t$ is needed, for which the tools developed in Sections~\ref{sec:R1} and \ref{sec:smallt} may be a starting point.

\section*{ACKNOWLEDGEMENTS}

We thank Yosef Rinott for many helpful discussions on this topic.  Bartroff's work was supported by grant DMS-0907241 from the National Science Foundation.

\end{document}